\documentclass[a4paper]{article}

\catcode`\é=13  \defé{\'e}
\catcode`\â=13  \defâ{\^a}
\catcode`\è=13  \defè{\`e}
\catcode`\à=13  \defà{\`a}
\catcode`\ù=13  \defù{\`u}
\catcode`\ç=13  \defç{\c c}
\catcode`\ê=13  \defê{\^e}
\catcode`\ô=13  \defô{\^o}
\catcode`\û=13  \defû{\^u}
\catcode`\î=13  \defî{\^\i}
\usepackage{amssymb}
\usepackage{amsthm}
\usepackage[english,francais]{babel}
\newtheorem{theorem}{Theorem}[section]

\newtheorem{e-proposition}[theorem]{Proposition}

\newtheorem{e-definition}[theorem]{Definition\rm}

\newtheorem{theoreme}{Th\'eor\`eme}[section]

\newtheorem{proposition}[theoreme]{Proposition}

\def\og{\leavevmode\raise.3ex\hbox{$\scriptscriptstyle\langle\!\langle$~}}
\def\fg{\leavevmode\raise.3ex\hbox{~$\!\scriptscriptstyle\,\rangle\!\rangle$}}

\begin{document}
\selectlanguage{francais}
\title{Variance limite d'une marche aléatoire réversible en milieu aléatoire sur~${\mathbb Z}$%
\footnote{A paraître aux Comptes Rendus de l'Académie des Sciences, Paris.}
\\ \normalsize\bf Limit of the  Variance of a Reversible Random Walk in  Random Medium on ${\mathbb Z}$.}
\author{Jérôme Depauw\footnote{Laboratoire de Mathématique et Physique Théorique, UMR CNRS 6083, Université Rabelais, Parc de Grandmont, 37200 Tours, France}, Jean-Marc Derrien\footnote{Laboratoire de Mathématiques de Brest,
 UMR CNRS 6205, Université de Bretagne Occidentale, 6, avenue Victor Le Gorgeu, CS 93837, F-29238 BREST Cedex 3
}}
\maketitle

\begin{abstract}
 Le théorème limite central pour la marche aléatoire sur un réseau aléatoire stationnaire de conductances a été étudié par de nombreux auteurs. En dimension 1, lorsque conductances et résistances sont inté\-gra\-bles, on peut montrer, pour presque tout environnement, la convergence
vers une loi gaussienne non dégénérée en suivant une méthode de
martingales introduite
dans ce contexte 
par S. Kozlov (1985).
Lorsque les résistances
ne sont pas intégrables,  Y. Derriennic et M. Lin 
 ont établi la convergence,
 cette
fois avec
 variance nulle, et  en probabilité relativement aux environnements
(communication personnelle). On montre ici, par une méthode particulièrement simple, que cette dernière convergence a lieu ponctuellement. Le problème analogue pour la diffusion continue est ensuite considéré. Enfin 
notre méthode nous permet de démontrer
une inégalité sur la moyenne quadratique d'une diffusion $(X_t)_t$, à temps $t$ fini.

\selectlanguage{english}
\noindent\begin{center}
{\bf Abstract}
\end{center}

The Central Limit Theorem for the random walk on a stationary random network of conductances
 has been studied by several authors. In one dimension, when conductances and resistances are integrable,
and following a  method of martingale introduced by S. Kozlov (1985), we
can prove the Quenched Central Limit Theorem.
In that case the variance of the limit law is not null. When resistances are not integrable, the Annealed
 Central Limit Theorem with null variance was established by Y. Derriennic and
 M. Lin (personal communication). 
 The quenched version of
this last theorem is proved here, by using a very simple method. The similar problem for the continuous diffusion is then considered.
 Finally our method allows us to prove an
 inequality for the quadratic mean of a diffusion $(X_t)_t$ at all time $t$.
\end{abstract}

\section*{Abridged English version}

We consider, on the ${\mathbb Z}$-network, a random stationary sequence of conductances, defined through a probability
space $(\Omega,{\cal A},\mu)$, an invertible    $\mu$-preserving   transformation $T$ which is also  ergodic, and a random variable $c>0$:
for a fixed environment $\omega\in\Omega$, the conductance of the edge
$[k,k+1]$ is $ c(T^k\omega)$.
Let $\bar c= c+c\circ T^{-1}$. 
We introduce the random walk $(X_n)_n$ on ${\mathbb Z}$ with Markov's operator $f\mapsto P_\omega f$ defined by
$$
P_\omega f(k)=\frac{1}{\bar c(T^k\omega)}\Bigl(c(T^{k-1}\omega)f(k-1)+c(T^k\omega)f(k+1)\Bigr),
$$ 
and with initial condition  $X_0=0$.

If  $c$ and $c^{-1}$ are $\mu$-integrable, we have a Quenched Central Limit Theorem: for almost all environment
$\omega\in\Omega$ the randon variable $\frac{X_n}{\sqrt{n}}$ converges in law, for $n\to\infty$, to the Gaussian law with mean $0$ and variance $\sigma^2=\displaystyle \Bigl[ \Bigl(\int_\Omega c\ d\mu\Bigr)\Bigl(\int_\Omega\frac{1}{c}\ d\mu\Bigr)\Bigr]^{-1}.
$ 
To prove this result we can
use a modification of the walk $(X_n)_n$, to obtain a martingale. This
method is explained in~\cite{Kozlov} in the multidimensional case. The above
expression of the variance $\sigma^2$,
only true in dimension 1,
appears in~\cite{DeMasi;Ferrari;Goldstein;Wick}, where the Annealed
 Central Limit Theorem
is proved.

When $c$ is integrable but not $c^{-1}$,  Y. Derriennic and M. Lin have proved, in a not published work, the Annealed
 Central Limit Theorem with null variance:
 $\displaystyle \lim_{n\rightarrow+\infty}n^{-1}{E_\omega(X_n^2)}=0$ in $\mu$-measure, where $E_\omega$ denotes the expectation relatively
 to the randomness of the walk, the environment being fixed.
We prove here the quenched version, i.e. the $\mu$-a.e. convergence.
More precisely, we prove the following theorem, which gives the value of  $\displaystyle \lim_{n\rightarrow+\infty}n^{-1}{E_\omega(X_n^2)}$, without
 condition on $c$ (other than $c>0$).
\begin{theorem}\label{eth1}
We have, for almost all environments $\omega$, 
$$ \lim_{n\rightarrow+\infty}\frac{E_\omega(X_n^2)}{n}=\displaystyle \Bigl[ \Bigl(\int_\Omega c\ d\mu\Bigr)\Bigl(\int_\Omega\frac{1}{c}\ d\mu\Bigr)\Bigr]^{-1},
$$
This limit being null if one of the integrals is $+\infty$ (or  both).
\end{theorem}
M. Biskup, T.M. Prescott in~\cite{Biskup;Prescott} and P. Mathieu in~\cite{Mathieu} have obtained a quenched functional
central limit theorem
in ${\mathbb Z}^d$, $d\geq 2$, for random walks in random media of uniformly bounded conductances
which form a family of independent identically distributed random variables.
The limiting brownian motion is always non-degenerate even if the resistances are not
integrable.
\par\noindent{\it Proof. ---}
Our method is particularly simple and does not use any martingale. On the other hand, it does not give the 
  Central Limit Theorem.
But our computation is still valid with infinite integrals.
Fix $\omega $ and denote by $P$ the operator denoted by $P_\omega$ before.  Consider  a function $f$, defined on ${\mathbb Z}$, satisfying
$
(P-I)f\equiv 1$, and
 $f(0)=0$. 
We can take for instance
$$f(m)=\left\{\matrix{
\displaystyle\sum_{\ell=0}^{m-1}\frac{1}{c(T^\ell\omega)}\sum_{k=1}^\ell \bar c(T^k\omega)&\mbox{ if }&m\geq 1,\cr
\displaystyle\sum_{\ell=1}^{-m}\frac{1}{c(T^{-\ell}\omega)}\sum_{k=0}^{\ell-1} \bar c(T^{-k}\omega)&\mbox{ if }&m\leq -1.}\right.
$$
It is easy to deduce from the pointwise ergodic theorem (See~\cite{Wiener}) that
\begin{equation}\label{eE2}
\lim_{m\rightarrow \pm\infty}\frac{f(m)}{m^2}=\frac{1}{2}\Bigl(\int_\Omega \frac{1}{ c}\ d\mu\Bigr)\Bigl(\int_\Omega \bar c \ d\mu\Bigr)
\qquad \mbox{p.s.}\ \omega.
\end{equation}
The point is that, since $c>0$, this convergence is still satisfied if one of these integrals is $+\infty$ (or  both).
Moreover, from the definition of $f$, we have $E_\omega(f(X_n))=n$, which can be rewritten as
$
\displaystyle E_\omega\Bigl(\frac{f(X_n)}{X_n^2}\times\frac{X_n^2}{n}\Bigr)=1
$. Finally,
 by considering separately points where $|X_n|$ is either small or large following~(\ref{eE2}), we obtain
$$
\frac{1}{2}\Bigl(\int_\Omega \bar c \ d\mu\Bigr)\Bigl(\int_\Omega \frac{1}{ c}\ d\mu\Bigr) E_\omega\Bigl(\frac{X_n^2}{n}\Bigr)\longrightarrow 1
$$
for $n$ tending to infinity, which prove Theorem~\ref{eth1}.\qed
\par We can state an analogue of Theorem~\ref{eth1} for continuous time and discrete space, and an analogue statement 
for both continuous time and space.
\par
Changing our framework, we  consider now non random environment.  The  continuous analogue of equation $(P-I)f\equiv 1$ allows us
to prove the  following very natural result.   
\begin{e-proposition}
Let $\sigma:{\mathbb R}\rightarrow{\mathbb R}_+^*$  be  function with a locally Lipschitz first derivative $\sigma'$. Suppose that 
$\sigma$ is bounded by a constant $\sigma_0$. Then the solution $(X_t)_t$
of the stochastic differential equation ${\rm d}X_t=\sigma(X_t)\ {\rm d}B_t+\sigma(X_t)\sigma'(X_t)\ {\rm d}t$ satisfies $ E(X^2_t)\leq \sigma_0^2t$
for any $t>0$.
\end{e-proposition}
The similar statement for the equation without drift ${\rm d}X_t=\sigma(X_t)\ {\rm d}B_t$  is well-known.
\selectlanguage{francais}
\section{Introduction}
Soit, sur le réseau ${\mathbb Z}$, une suite aléatoire stationnaire de conductances, définie à l'aide d'un espace probabilisé
$(\Omega,{\cal A},\mu)$, d'une transformation inversible et ergodique $T$ préservant la mesure $\mu$, et d'une variable aléatoire $c>0$:
pour un environnement $\omega\in\Omega$ fixé, la conductance de l'arête 
$[k,k+1]$ est $ c(T^k\omega)$.
On pose $\bar c=c+T^{-1}c$. 
Soit alors
 la marche aléatoire $(X_n)_n$ sur ${\mathbb Z}$ d'opérateur de Markov $f\mapsto P_\omega f$ défini par
$$
P_\omega f(k)=\frac{1}{\bar c(T^k\omega)}\Bigl(c(T^{k-1}\omega)f(k-1)+c(T^k\omega)f(k+1)\Bigr),
$$ 
et de condition initiale $X_0=0$.

Lorsque les fonctions $c$ et $c^{-1}$ sont intégrables, on a un théorème limite central pour presque tout environnement
$\omega\in\Omega$: la variable aléatoire $\frac{X_n}{\sqrt{n}}$ tend en loi, pour $n$ tendant vers l'infini, vers la loi de Gauss de moyenne nulle et de variance $\sigma^2=\displaystyle \Bigl[ \Bigl(\int_\Omega c\ d\mu\Bigr)\Bigl(\int_\Omega\frac{1}{c}\ d\mu\Bigr)\Bigr]^{-1}$. 
On peut démontrer ce résultat en suivant la méthode de martingales développée dans~\cite{Kozlov}  dans un
cadre multidimensionnel.
 Enfin l'expression
de $\sigma^2$ ci-dessus, typique de la dimension 1, apparaît
dans~\cite{DeMasi;Ferrari;Goldstein;Wick}, où la convergence en
moyenne relativement aux environnements est démontrée.

Une étude approfondie de cette même méthode, non publiée, a permis ré\-cem\-ment  à  Y. Derriennic et M. Lin de démontrer que, 
lorsque les conductances  $c$ sont  intégrables mais pas les résistances
$c^{-1}$, 
 on a la convergence  
 $\displaystyle \lim_{n\rightarrow+\infty}n^{-1}{E_\omega(X_n^2)}=0$ en probabilité pour $\mu$, la notation $E_\omega$
désignant l'es\-pé\-ran\-ce relativement à l'aléa de la chaîne, à environnement fixé.
On montre ici le résultat analogue, pour la convergence presque sûre.
Plus précisément, on montre le théorème suivant, qui donne la valeur de $\displaystyle \lim_{n\rightarrow+\infty}n^{-1}{E_\omega(X_n^2)}$, sans condition sur $c$ (autre que $c>0$).
\begin{theoreme}\label{th1}
On a, pour presque tout environnement $\omega$, 
$$ \lim_{n\rightarrow+\infty}\frac{E_\omega(X_n^2)}{n}=\displaystyle \Bigl[ \Bigl(\int_\Omega c\ d\mu\Bigr)\Bigl(\int_\Omega\frac{1}{c}\ d\mu\Bigr)\Bigr]^{-1},
$$
cette limite étant nulle dès que l'une des deux intégrales diverge (ou les deux).
\end{theoreme}
M. Biskup, T.M. Prescott dans~\cite{Biskup;Prescott} et P. Mathieu dans~\cite{Mathieu} ont 
obtenu un principe d'invariance ponctuel dans ${\mathbb Z}^d$, $d\geq 2$, pour des conductances
bornées, indépendantes et équidistribuées.
Le mouvement brownien limite est non dé\-gé\-né\-ré, même si les résistances ne sont pas intégrables.

\section{Démonstration}
Notre méthode, particulièrement simple, ne repose pas sur un argument de type martingale. Elle ne donne pas le 
 théorème limite central.
Par contre notre calcul de la
variance reste valide dans le cas dégénéré.
Soit $\omega $  fixé. On note $P$ l'opérateur noté $P_\omega$ ci-dessus. On considère  une fonction $f$, définie sur ${\mathbb Z}$, vérifiant
$
(P-I)f\equiv 1$, et
 $f(0)=0$.
On a donc $E_\omega(f(X_n))-E_\omega(f(X_{n-1}))=1$ et $E_\omega(f(X_0))=0$, soit encore 
\begin{equation}\label{E1}
E_\omega(f(X_n))=n.
\end{equation}
D'autre part, en notant $\tau$ la translation des fonctions, définie par $(\tau f)(k)=f(k+1)$, et en considérant
 $\partial =I-\tau$ et $\partial^*=I-\tau^{-1}$,  on a
 $\displaystyle I-P=\frac{1}{\bar c}\partial^*c\partial$. Il en découle aisément l'expression d'une solution $f\geq 0$ explicite: $f(0)=f(1)=0$, et
$$f(m)=\left\{\matrix{
\displaystyle\sum_{\ell=0}^{m-1}\frac{1}{c(T^\ell\omega)}\sum_{k=1}^\ell \bar c(T^k\omega)&\mbox{ if }&m\geq 1,\cr
\displaystyle\sum_{\ell=1}^{-m}\frac{1}{c(T^{-\ell}\omega)}\sum_{k=0}^{\ell-1} \bar c(T^{-k}\omega)&\mbox{ if }&m\leq -1.}\right.
$$
Il est facile de voir que le théorème ergodique ponctuel (voir~\cite{Wiener}), appliqué successivement aux deux sommations en $k$ et $\ell$, donne alors
\begin{equation}\label{E2}
\lim_{m\rightarrow \pm\infty}\frac{f(m)}{m^2}=\frac{1}{2}\Bigl(\int_\Omega \frac{1}{ c}\ d\mu\Bigr)\Bigl(\int_\Omega \bar c \ d\mu\Bigr)
\qquad \mbox{p.s.}\ \omega
\end{equation}
Le point est que, comme $c>0$, cette convergence est aussi vérifiée si l'une des deux intégrales diverge (ou les deux), la limite étant alors $+\infty$.
Notons $v_\infty$ l'inverse multiplicatif, éventuellement nul, de cette limite.
Finalement de l'égalité~(\ref{E1}),
on tire 
\newpage
\begin{eqnarray*}
E_\omega\Bigl(\frac{X_n^2}{n}\Bigr)-v_\infty&=&
E_\omega\Bigl(\Bigl(\frac{X_n^2}{f(X_n)}-v_\infty\Bigr)\frac{f(X_n)}{n}{\mathbf 1}_{|X_n|> M}\Bigr)+\cr
&&+\frac{1}{n}E_\omega\Bigl(\Bigl(X_n^2-v_\infty f(X_n)\Bigr){\mathbf 1}_{|X_n|\leq M}\Bigr).
\end{eqnarray*}
Pour $n$ tendant vers l'infini, le  dernier terme s'annule. La
démonstration du théorème~\ref{th1} s'achève en laissant ensuite
tendre  $M$ vers l'infini, et en appliquant les égalités~(\ref{E1}) et~(\ref{E2}).\qed
\section{Analogues continus}
Le modèle précédent de la marche aléatoire sur un réseau de conductances aléatoires a naturellement deux analogues en paramètres continus. 
\begin{itemize}
\item Le premier analogue consiste à prendre
le temps continu mais  l'espace toujours  discret. On considère
 alors, en général,
 le processus de Markov $(X_t)_{t\in{\mathbb R}}$ sur ${\mathbb Z}$
 de générateur infinitésimal
\begin{equation}\label{Lo}
L_\omega f(k)=c(T^{k-1}\omega)f(k-1)+c(T^k\omega)f(k+1)-\bar c(T^k\omega)f(k),
\end{equation}
soit encore $L_\omega f=-\partial^*(c\partial f)$. Par une démonstration  très similaire à celle du théorème~\ref{th1}, consistant à considérer une solution $f$ de l'équation
$L_\omega f\equiv 1$, on obtient comme limite de la variance 
$\displaystyle \lim_{t\rightarrow+\infty}t^{-1}E_\omega(X_t^2)=2\Bigl[\displaystyle\int_\Omega c^{-1}\ d\mu\Bigr]^{-1},
$
cette limite étant nulle dès que l'intégrale diverge (le théorème limite central dans le cas où $c^{-1}$ est intégrable est dû à K. Kawazu et H. Kesten; 
voir~\cite{Kawazu;Kesten}).
\item   Le second analogue consiste à prendre à la fois le temps et l'espace continus. Ce modèle étant 
 formellement moins similaire,
et plus intuitif, nous
nous y attardons un peu plus dans le paragraphe suivant.
\end{itemize}
\section{Diffusion en milieu aléatoire stationnaire}
La version continue, en temps et en espace, du travail précédent,
 consiste à se donner un espace probabilisé $(\Omega,{\cal A},\mu)$,
 muni d'un flot $(T_x)_{x\in{\mathbb R}}$ ergodique
préservant la probabilité $\mu$. Précisément, on suppose
 que l'application $(\omega,x)\mapsto T_x\omega$
est mesurable et vérifie:
\begin{itemize}
\item $T_{x+y}=T_xT_y$ et $T_0\omega=\omega$;
\item si, pour tout $x\in{\mathbb R}$, on a $T_x F=F$ modulo $\mu$, alors $\mu(F)=0$ ou $1$; 
\item $\mu(T_x F)=\mu(F)$.
\end{itemize}
On se donne aussi deux variables aléatoires $\lambda$, $\gamma>0$, et on
 étudie, à $\omega$ fixé, 
 le processus de générateur infinitésimal $f\mapsto L_\omega f$ défini par
\begin{equation}\label{gen}
L_\omega f(x)=\frac{1}{2\gamma (T_x\omega)}\frac{{\rm d}}{{\rm d} x}\Bigl(\lambda(T_x\omega)\frac{{\rm d} f}{{\rm d} x}\Bigr),
\end{equation} 
et de condition initiale $X_0=0$.
L'interprétation physique est celle 
d'un conducteur thermique linéaire,  de longueur infinie,  tel que,  pour un environnement $\omega\in\Omega$,
la conductivité et la capacité thermiques soient données respectivement par les fonctions
$x\mapsto\lambda(T_x\omega)$ et $x\mapsto\gamma(T_x\omega)$.
Il est usuel d'écrire le problème sous forme d'une équation différentielle stochastique
 ${\rm d}X_t=\sigma_\omega(X_t)\ {\rm d}B_t+b_\omega(X_t)\ {\rm d}t$.
Ici le  coefficient de diffusion est donné (toujours pour un environement $\omega$ fixé) par $\sigma_\omega^2(x)=\displaystyle \lambda(T_x\omega)\gamma(T_x\omega)^{-1}$, et la dérive par $b_\omega(x)=
\displaystyle (2\gamma(T_x\omega))^{-1}\frac{d}{dx} \lambda(T_x\omega)$.
Le théorème analogue au résultat du paragraphe précédent prend la forme suivante.
\begin{theoreme}\label{th2}
Supposons que, pour presque tout $\omega\in\Omega$, les fonctions
$\sigma_\omega^2$ et $b_\omega$ sont localement lipschitziennes. Alors,
pour presque tout environnement $\omega$, la solution $(X_t)_t$ de l'équation différentielle stochastique ci-dessus vérifie
$$ \lim_{t\rightarrow+\infty}\frac{E_\omega(X_t^2)}{t}=\displaystyle \Bigl[\Bigl(\int_\Omega \gamma\ d\mu\Bigr)\Bigl(\int_\Omega\frac{1}{\lambda}\ d\mu\Bigr)\Bigr]^{-1},
$$
cette limite étant nulle dès que l'une des deux intégrales diverge (ou les deux).
\end{theoreme}
Là encore, dans le cas où les fonctions $\gamma$ et $\lambda^{-1}$ sont intégrables sur $\Omega$,
 il est connu que
le théorème limite central est valide (voir notamment G.D. Papanicolaou et S.R.S. Varadhan~\cite{Papanicolaou;Varadhan} pour le cas elliptique).

\par\noindent{\it Démonstration du théorème~\ref{th2}.---} La démonstration est très similaire à celle 
du cas de la marche aléatoire. L'environnement $\omega\in\Omega$ étant fixé,
on considère la fonction $f\geq 0$ définie sur ${\mathbb R}$ par 
$$f(x)=\left\{\matrix{\displaystyle
\int_{v=0}^x\frac{1}{\lambda(T_v\omega)}\int_{u=0}^v2\gamma(T_u\omega)\ {\rm d}u\
{\rm d}v&\mbox{ si }&x\geq 0\cr
\displaystyle\int_{v=x}^0\frac{1}{\lambda(T_v\omega)}\int_{u=v}^02\gamma(T_u\omega)\ {\rm d}u\
{\rm d}v &\mbox{ si }&x\leq 0.}\right.
$$
Le théorème ergodique ponctuel, appliqué successivement aux deux intégrales en $u$ et $v$, donne 
$$
\lim_{x\rightarrow \pm\infty}\frac{f(x)}{x^2}=\Bigl(\int_\Omega \frac{1}{ \lambda}\ d\mu\Bigr)\Bigl(\int_\Omega \gamma \ d\mu\Bigr)
\qquad \mbox{p.s.}\ \omega,
$$
et comme les fonctions sont positives, la convergence est aussi vérifiée si l'une des deux intégrales diverge (ou les deux), la limite étant alors $+\infty$.
Comme d'autre part, d'après les hypothèses du théorème, la fonction $f$ est de classe
${\cal C}^2$, le processus $Y$  défini par $Y_t=f(X_t)$ vérifie une équation différentielle stochastique obtenue grâce à la formule de Itô.
Le calcul, que nous ne développons pas, donne un coefficient de dérive constant: $
{\rm d}Y_t=a_\omega(X_t)\ {\rm d}B_t+\ {\rm d}t.
$
On a donc $E_\omega(f(X_t))=t$, et la démonstration se termine comme celle du théorème~\ref{th1}.\qed
\section{Majoration et minoration de la moyenne quadratique à temps $t$ fini}
Il est clair, d'après l'expression de la variance limite de la marche aléatoire donnée par le théorème~\ref{th1},
que la variance de $X_n$ n'est pas une fonction croissante de la conductance $c$. Par contre, en temps continu, la question de la monotonie
 de $E_\omega(X_t^2)$ par rapport à chacun des coefficients
 de capacité et de conductivité se pose. 
 Considérons ici la dépendance par rapport à la conductivité.
 Dans le générateur infinitésimal~(\ref{gen}),
 prenons $\gamma\equiv 1$, et
 notons $\sigma^2$  la fonction, toujours supposée strictement positive, donnant le coefficient de diffusion, soit $\sigma^2(\omega)
 =\lambda(\omega)$. A défaut de véritable loi de monotonie, on a le résultat suivant.
 \begin{proposition}\label{th3}
Supposons que pour presque tout $\omega\in\Omega$, la fonction
$x\mapsto  \sigma^2(T_x\omega)$ est dérivable, de dérivée localement
lipschitzienne. On suppose éga\-le\-ment qu'il existe une constante $\sigma_0^2>0$ telle que 
p.s.~$\omega$,  on a $ \sigma^2(\omega)\leq
 \sigma_0^2$. Alors,
pour presque tout environnement $\omega$ on a, 
pour tout $t$, $E_\omega(X_t^2)\leq \sigma_0^2t.
$
\end{proposition}
Comme cela se voit dans la démonstration ci-dessous, cette proposition n'a rien à voir avec l'environnement aléatoire. Ce cadre
 n'est gardé ici que pour éviter
d'introduire de nouvelles notations. Ce résultat est très naturel, et sa démonstration est une très simple application de la méthode utilisée dans la
 démonstration ci-dessus (l'inégalité équivalente pour $\gamma$ variable mais $\lambda\equiv 1$, soit une dérive $b=0$, est classique). 
\par\noindent{\it Démonstration. ---} La fonction $f$ utilisée dans la démonstration du théorème~\ref{th2} vérifie maintenant 
 $f(x)\geq \sigma_0^{-2}x^2$. Comme d'autre part  $E_\omega(f(X_t))=t$, la proposition~\ref{th3}
est démontrée.\qed
\par
Lorsque l'on suppose $ \sigma^2(\omega)\geq
 \sigma_0^2$, on a de même la minoration ${E_\omega(X_t^2)}\geq \sigma_0^2t$. 

\bibliographystyle{plain}

\end{document}